# Scharlemann's manifold is standard

By Selman Akbulut*

*Dedicated to Robion Kirby on the occasion of his 60th birthday*

## Abstract

In his 1974 thesis, Martin Scharlemann constructed a fake homotopy equivalence from a closed smooth manifold $f : Q \to S^3 \times S^1 \# S^2 \times S^2$, and asked the question whether or not the manifold $Q$ itself is diffeomorphic to $S^3 \times S^1 \# S^2 \times S^2$. Here we answer this question affirmatively.

In [Sc] Scharlemann showed that if $\Sigma^3$ is the Poincaré homology 3-sphere, by surgering the 4-manifold $\Sigma \times S^1$, along a loop in $\Sigma \times 1 \subset \Sigma \times S^1$ normally generating the fundamental group of $\Sigma$, one obtains a closed smooth manifold $Q$ and homotopy equivalence:

$$f : Q \longrightarrow S^3 \times S^1 \# S^2 \times S^2$$

which is not homotopic to a diffeomorphism (actually by taking $\Sigma$ to be any Rohlin invariant 1 homology sphere one gets the same result). He then posed the question whether $Q$ is a standard copy of $S^3 \times S^1 \# S^2 \times S^2$ (i.e. whether $f$ is a fake self-homotopy equivalence) or $Q$ itself is a fake copy of $S^3 \times S^1 \# S^2 \times S^2$.

This question has stimulated much research during the past twenty years resulting in some partial answers. For example, in [FP] it was shown that $Q$ is stably standard, in [Sa] it was proven that it is obtained by surgering a knotted $S^2 \subset S^2 \times S^2$, and in [A4] it was shown that it is obtained by "Gluck twisting" $S^3 \times S^1 \# S^2 \times S^2$ along an imbedded 2-sphere. Also by a similar construction one can obtain a fake homotopy equivalence:

$$g : P \longrightarrow S^3 \tilde{\times} S^1 \# S^2 \times S^2$$

where $S^3 \tilde{\times} S^1$ is the nonorientable $S^3$ bundle over $S^1$. But in this case it turns out that $P$ is not standard, i.e. $P$ is a fake copy of $S^3 \tilde{\times} S^1 \# S^2 \times S^2$ which is also obtained by Gluck twisting a 2-sphere in $S^3 \tilde{\times} S^1 \# S^2 \times S^2$ ([A2], [A3], [A4]).

---

*Partially supported by NSF grant DMS-9626204 and the research at MSRI is supported by in part by NSF grant DMS-9022140.

This is the only example known to the author where one can make a smooth 4-manifold fake, by Gluck twisting an imbedded 2-sphere. Since Gluck twisting operation preserves gauge theoretical invariants of oriented 4-manifolds, it is hard to find oriented such examples (this is because the manifold obtained by this operation is stably diffeomorphic to the standard one, under connected summing operation with either $\pm \mathbb{CP}^2$).

Even though by using surgery techniques it was observed by R. Lee that $S^3 \times S^1 \# S^2 \times S^2$ does in fact admit a fake self-homotopy equivalence (see [CS1], p. 515), Scharlemann's manifold remained a source of hope for topologist as a possible way to establish the existence of a fake $S^3 \times S^1$. Here we prove:

THEOREM.  *Q is diffeomorphic to $S^3 \times S^1 \# S^2 \times S^2$.*

Note that by Rohlin's theorem the Poincaré homology sphere $\Sigma$ cannot imbed into either $S^3 \times S^1$ or $S^2 \times S^2$. However, an immediate corollary of the Theorem is that $\Sigma$ embeds into $S^3 \times S^1 \# S^2 \times S^2$

In this paper we use the convention of [A1], denoting a 1-handle with a "circle with dot."

*Constructions and the proof.* We first write $\Sigma \times S^1 = \Sigma \times I_+ \smile \Sigma \times I_-$, where $I_\pm \approx I = [0, 1]$ are closed intervals and the union is taken along the boundaries (i.e. along $\Sigma \sqcup -\Sigma$). Let $N$ be the manifold obtained by surgering a loop in $\Sigma \times I_-$ normally generating the fundamental group of $\Sigma$. Hence $Q = \Sigma \times I_+ \smile_\partial N$. We will now construct a handlebody of $Q$ (Figs. 11, 12 and 14): Figure 1 is the $B^4$ with a 2-handle attached to the left-handed trefoil knot $K$ with $-1$ framing. It is well known that the boundary of this manifold is the Poincaré homology sphere $\Sigma$.

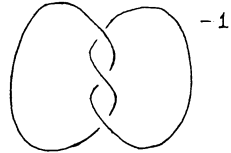

Figure 1.

*Handlebody of $\Sigma \times I$.* We now want to visualize $\Sigma \times I_-$: This manifold is obtained by removing the tubular neighborhood of the trefoil knot $K$ from $S^3$ and crossing with $I$ and attaching a 2-handle to the boundary of this 4-manifold with $-1$ framing push-off of the trefoil knot, which we call $\alpha$.

Since 3 and 4 handles are always attached in the standard way, we do not need to visualize them, we only need to describe 4-dimensional handlebodies by indicating their 1 and 2 handles (together with the knowledge of whether there are 3 and 4 handles).

Up to attaching a 3 handle, $\Sigma \times I_-$ is obtained by removing the tubular neighborhood of a properly imbedded arc (with trefoil knot tied on it) from $B^3$, and crossing with $I$, and attaching a 2-handle to the boundary of this 4-manifold along the "$-1$ framing push-off" of the trefoil knot $\alpha$, as indicated in Figure 2. Clearly, this is obtained by removing the "usual" slice disc from $B^4$, which the trefoil knot connected summed with its mirror image $K\#(-K)$ bounds, and attaching the 2-handle $\alpha$, as indicated in Figure 3. The dot on the knot $K\#(-K)$ in Figure 3 indicates that the tubular neighborhood of the slice disc which it bounds is removed from $B^4$. We will refer this as *slice 1-handle*. This notation was discussed in [AK1], for example either of the pictures in Figure 4 describes the handlebody (two 1 and one 2-handles) of the $B^4$ with this slice disc removed (canceling one of the 1-handles of Figure 4 by the 2-handle gives the slice 1-handle of Figure 3). Now on the boundary of the handlebody of Figure 4, framed knot $\alpha$ sits as indicated in Figure 5. This can be checked by keeping track of $\alpha$ while proceeding from Figure 3 to Figure 4 (at this stage reader should disregard linking $+1$ and $0$ framed handles added to Figure 5 as well as the "arrow"; these will be explained in the forthcoming steps).

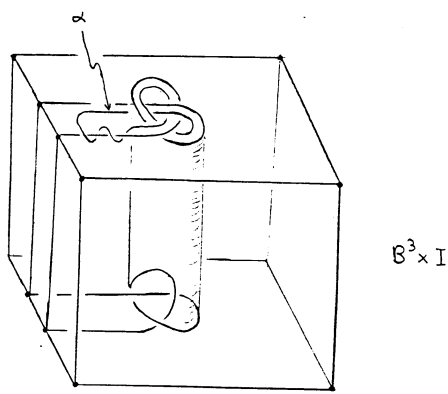

Figure 2.

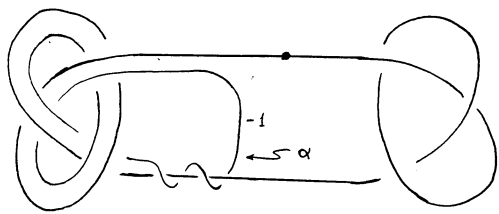

Figure 3.

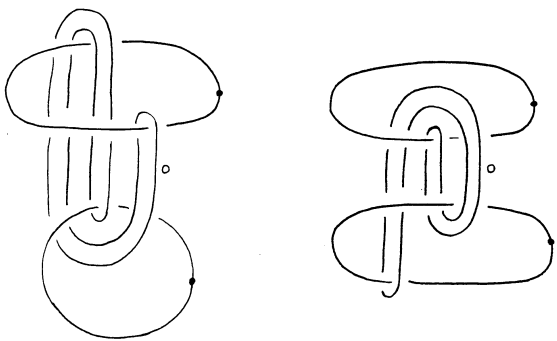

Figure 4.

*Handlebody of $N$.* Recall that surgering a loop from 4-manifold (i.e. cutting out a $B^3 \times S^1$ and gluing in $S^2 \times B^2$) corresponds to attaching a pair of 2-handles; one to this loop with any framing $k$, and the other to the linking circle of this loop with 0-framing. Clearly the diffeomorphisim type of the surgered manifold depends only on the parity of $k$ ($k$ can be changed by 2, by sliding it over the 0-framed handle). For example changing parity could change a spin manifold to nonspin manifold.

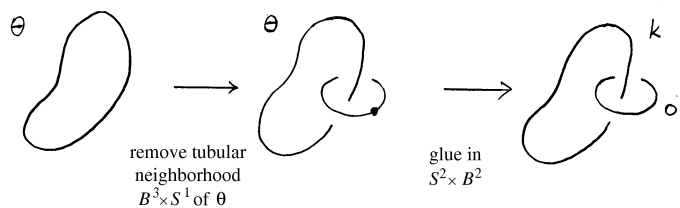

Now we surger $\Sigma \times I_-$ by attaching a pair of linked $+1$ and $0$ framed 2-handles as shown in Figure 5 (i.e. two circles linked as the Hopf link in the figure). By sliding the 2-handle $\alpha$ over the 1-handle (as indicated by the arrow) we obtain Figure 6. Also going from Figure 5 to Figure 6, by sliding 0-framed handle of the surgery (one of the Hopf link circles) over the $+1$ framed handle, we turn it into a $-1$ framed handle (linking 1-handle). We then slide this $-1$ framed 2-handle over the 0-framed two handle connecting the two 1-handles, to get Figure 7. To obtain Figure 8 we simply slide $\alpha$ over the linking $-1$ framed handle (this unlinks $\alpha$ from the $-1$ framed handle, and decreases the linking of $\alpha$ with the 1-handle).

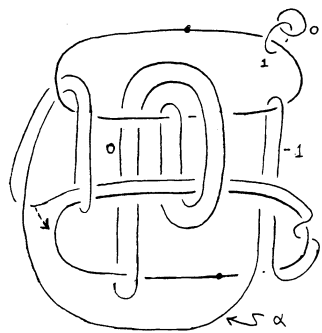 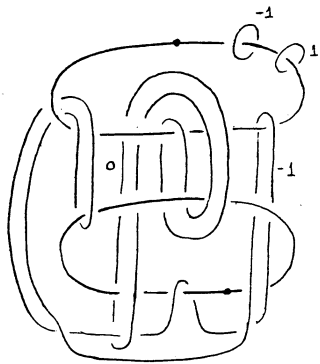

Figure 5.                    Figure 6.

Canceling the two 1-handles with the $-1$ and $+1$ framed 2-handles (which trivially link the 1-handles geometrically once) has an affect of twisting all the framed links going through the 1-handles by $-1$ and $1$ twists, respectively. Hence Figure 8 consists of basically pair of 2-handles attached to $B^4$. In Figure 8 we use the convention that the framed knots of the two 2-handles of Figure 8, going through the loops labeled by circled $\pm 1$, get twisted by $\pm 1$.

*Upside down handlebody of $N$.* We need to visualize $N$ as a handlebody "turned upside down," i.e. as handles attached to $\partial N = \Sigma \sqcup -\Sigma$ (more precisely to $\partial N \times I$).

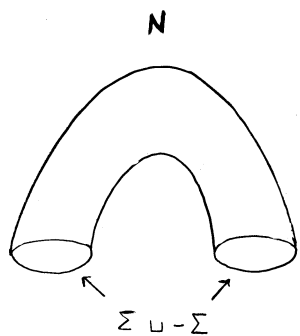

We claim Figure 9 describes this handlebody. Figure 9 is obtained by attaching one 1-handle and a pair of 0 framed 2-handles to $\Sigma \sqcup -\Sigma$. In Figure 9 $\Sigma \sqcup -\Sigma$ is drawn as a pair of trefoil knots with $\pm 1$ framings and a 3-handle (not drawn as usual). To check this we will turn this handlebody (we just described) "upside down," and see that we are getting Figure 8. To do this we simply attach 2-handles to the "dual loops" $\lambda$ and $\tau$ of the 2-handles of Figure 9 (these are small linking circles of the pair of 2-handles as indicated in the figure) to $B^4$.

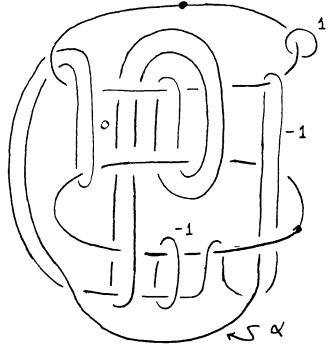 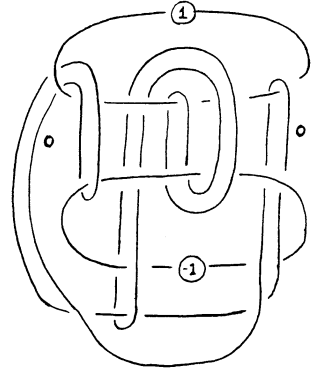

Figure 7.  Figure 8.

This is done by following the 0 framed loops $\lambda$ and $\tau$ through a diffeomorphism from the boundary of the handlebody of Figure 9 to $S^3$ (actually it suffices to get a diffeomorphism to $S^1 \times S^2$, since the upside down 1-handle of Figure 9 becomes a 3-handle turning $S^1 \times S^2$ to $S^3$ ), and then attaching 2-handles to $B^4$ along these framed loops $\lambda$ and $\tau$.

We construct a diffeomorphism from the boundary of Figure 9 to $S^3$ by canceling the pair of the two linking 0 framed 2-handles of the figure (and of course by canceling $S^1 \times S^2$ to $S^3$ by the 3-handle). This gives pair of disjoint unknotted $\pm 1$ framed circles in $S^3$ as in Figure 10, and by tracing through $\lambda$ and $\tau$, we see that they link these $\pm 1$ framed circles as in Figure 10. By blowing down these circles we get exactly Figure 8!

*Handlebody of $Q = \Sigma \times I \smile N$, and the Gluck twist.* We now attach the "upside down" handlebody of $N$ to $\Sigma \times I$ to obtain $Q$. This gives Figure 11. We now make a key observation: The attaching framed knot of the handle $\alpha$ is isotopic to the *trivially linking circle* of the slice 1-handle $K\#-K$, as indicated in Figure 12. This can be seen by applying the diffeomorphism to the boundary of either Figures 11 or 12, as described in the last paragraph, i.e. remove the "dot" from the slice 1-handle $K\#-K$, then cancel the pair of the 0 framed 2-handles. This turns the slice knot $K\#-K$ to an unknot and turns $\alpha$ to the desired (trivially linking) circle as in Figure 13. This means that the slice knot $K\#-K$ is unknotted on the boundary of $S^3 \times S^1 \# S^2 \times S^2 - \text{int}(B^4)$, and the slice disc it bounds along with the trivial disc it bounds in the 4-handle $B^4$ gives an imbedding:

$$S^2 \hookrightarrow S^3 \times S^1 \# S^2 \times S^2.$$

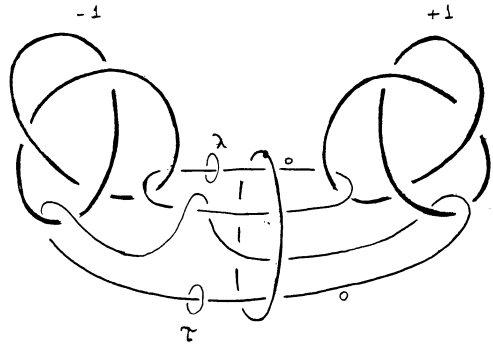

Figure 9.

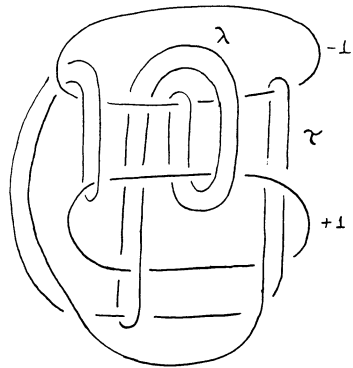

Figure 10.

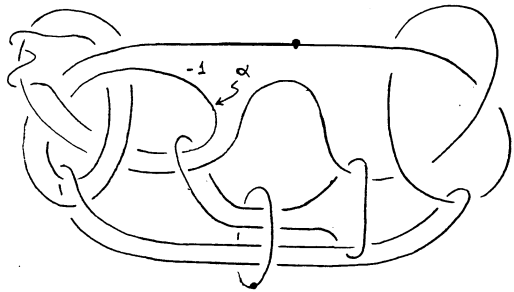

Figure 11.

504    SELMAN AKBULUT

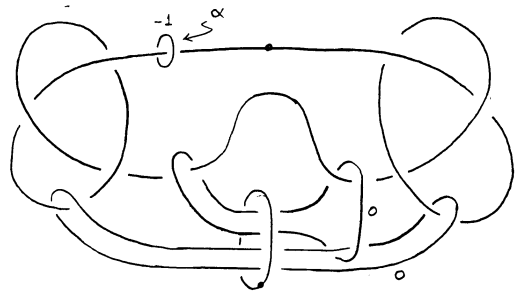

Figure 12.

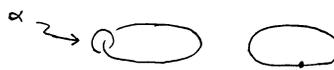

Figure 13.

Furthermore, the trivially linking $-1$ framed 2-handle $\alpha$ (of Figure 12) means that this imbedded $S^2$ has been Glucked (that is $S^3 \times S^1 \# S^2 \times S^2$ has been *Gluck twisted* along this 2-sphere), i.e. the tubular neighborhood $S^2 \times D^2$ of this 2-sphere has been removed and glued back by the nontrivial diffeomorphism of $S^2 \times S^1$ (see [A3],[A4]). One well-known fact about the Gluck construction is that the parity of framing of the "trivially linked circle" can be changed (e.g. the sign of the framing can be changed) without changing the diffeomorphism type of the manifold. For a quick reminder of a proof, we simply refer the reader to the following picture.

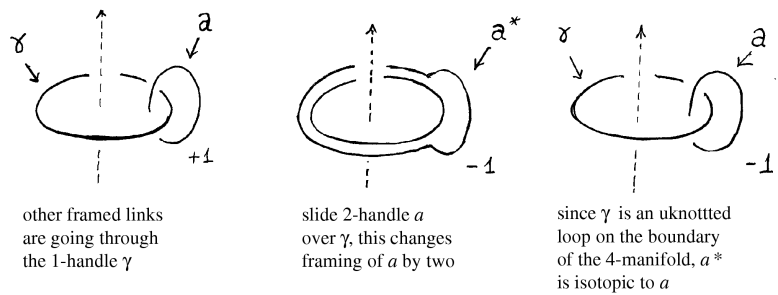

other framed links are going through the 1-handle $\gamma$

slide 2-handle $a$ over $\gamma$, this changes framing of $a$ by two

since $\gamma$ is an uknottted loop on the boundary of the 4-manifold, $a*$ is isotopic to $a$

This observation about Gluck construction was previously used in [A4] to establish the results about the Gluck twisting mentioned in the introduction. We redraw Figure 12 as Figure 14, by using the same steps as going from Figure 3 to Figure 4, to get a more concrete handlebody picture of $Q$.

*Canceling excess* 1-*handles of Q.* By sliding the $-1$ framed handle over the 0 framed handle as indicated in Figure 14, we obtain Figure 15. By sliding the 0 framed handle over the $-1$ framed handle as indicated in Figure 15, and then sliding the $-1$ framed handle over the 0 framed 2-handle (which connects the two 1-handles), we obtain Figure 16.

Now we get a pleasant surprise!: The two 1-handles of Figure 16 have been canceled by $\pm 1$ framed 2-handles, since they link the 1-handles geometrically once (to emphasize the canceling 2-handles, in Figure 16 arrows are drawn on their corresponding framed links). Now we can cancel these two 1- and 2-handle pairs, and be left with only one 1-handle and two 2-handles. We can now draw the picture of this simpler handlebody. Instead we will draw the upside down picture of this handlebody.

*Turning Q upside down.* To do this, as before, we simply take the dual 0 framed circles $\gamma$ and $\delta$ of the remaining (uncancelled) 2-handles of Figure 16, and trace them via a diffeomorphism of the boundary of Figure 16 to $\partial(B^3 \times S^1)$ and attach 2-handles to $B^3 \times S^1$ along these framed loops $\gamma$ and $\delta$ on the boundary. Now we proceed.

By sliding the $-1$ framed handle over the 0 framed handle (connecting the two 1-handles) as indicated in Figure 16, we obtain Figure 17. As shown in the figure, as a result of this move the position and framing of the loop $\gamma$ changed (e.g. its 0 framing has changed to $+1$ framing). Now by doing the reverse of the move we did going from Figure 14 to Figure 16 (namely slide $+1$ framed handle over the $-1$ framed handle, then slide $-1$ framed handle over the other 0 framed handle), we obtain Figure 18. Again, as a result this time the position and framing of the loop $\delta$ changed as shown in the figure (and its framing changed from 0 to $+1$).

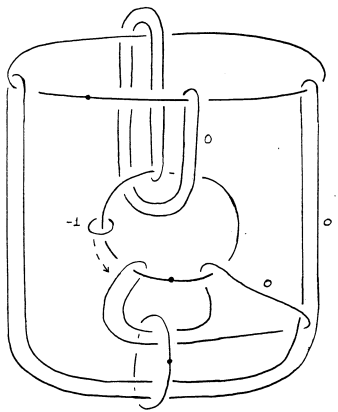
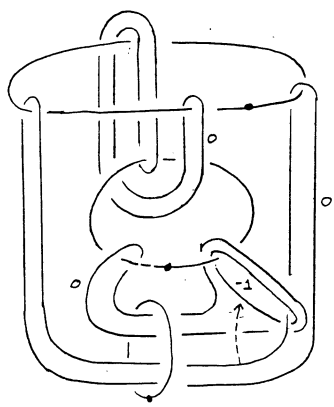

Figure 14.                                        Figure 15.

By applying Figure 18 the reverse of the move Figure 16 → Figure 17 we obtain Figure 19 and Figure 20. The reason we did this in two steps is to be careful about keeping track of the changing framing and the position of the loop $\gamma$ (now its framing is changed back to 0 and it goes back to its old position). Now we change the $-1$ framing of the "trivially linking circle" to the 1-handle to $+1$ framing (see the discussion in the section about Gluck twisting above).

We now apply the boundary diffeomorphism Figure 12 → Figure 13, i.e. we remove the dots from the 1-handles, by the aid of the connecting 0 framed handle turn them into $K \# - K$ , and cancel the two linking 0 framed 2-handles. The result is Figure 21. By sliding the two strands of the $+1$ framed handle over the 1-handle as indicated in Figure 21 we arrive to Figure 22.

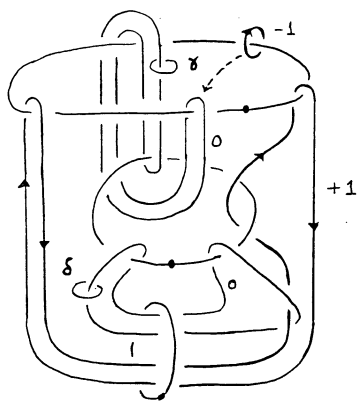
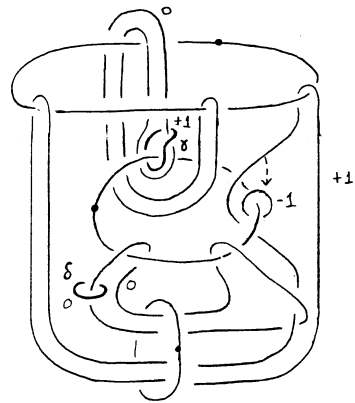

Figure 16.  Figure 17.

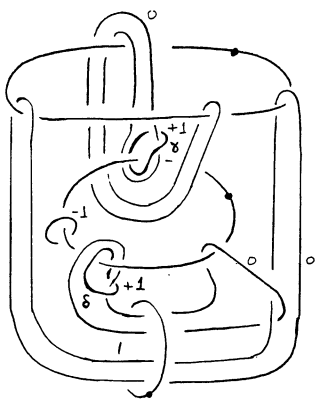
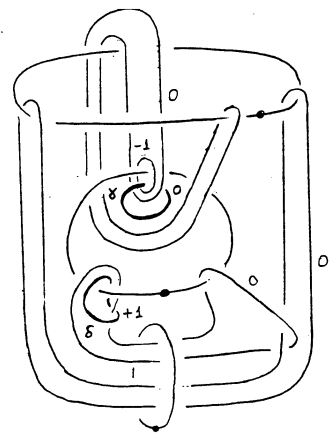

Figure 18.  Figure 19.

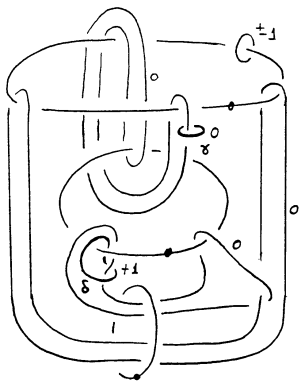 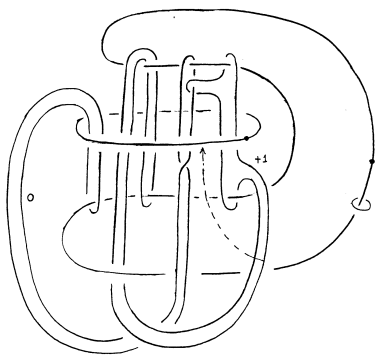

Figure 20.                    Figure 21.

We now repeat the by now standard move to Figure 22: By sliding the little +1 framed handle over the 0 framed handle we get Figure 23, then by sliding the other +1 framed handle over this +1 framed handle gives Figure 24. Notice that this move decreases linkings. By repeating this move once more, we arrive to Figure 25. A small isotopy of Figure 25 gives Figure 26. Then by sliding the 0 framed handle over the −1 framed handle we get Figure 27 which can be drawn as in Figure 28. By sliding the two strands of −1 framed handle over the +1 framed handle (as indicated in Figure 28) gives Figure 29.

By sliding the −1 framed handle over the −3 framed handle gives us Figure 30. Again by sliding the −1 framed handle three times over the 0 framed handle (as indicated in the figure) and letting the −3 framed handle slide over it we can completely undo the linking of the −3 framed handle and the 1-handle, by changing its framing to 0. This gives Figure 31. By canceling the obvious 1 and 2 handle pairs, we obtain Figure 32, and by undoing the knotted −4 handle by the 0 framed 2-handle we get Figure 33 which is $B^3 \times S^1 \# S^2 \times S^2$. Along with the 3 and 4 handles which we have been carrying along, this manifold actually is $S^3 \times S^1 \# S^2 \times S^2$. □

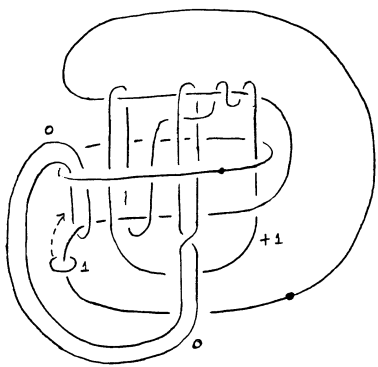 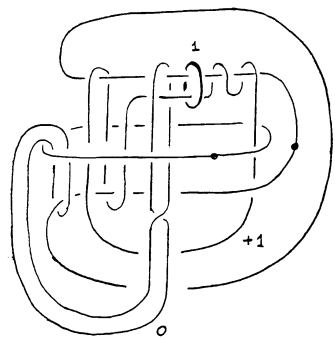

Figure 22.                    Figure 23.

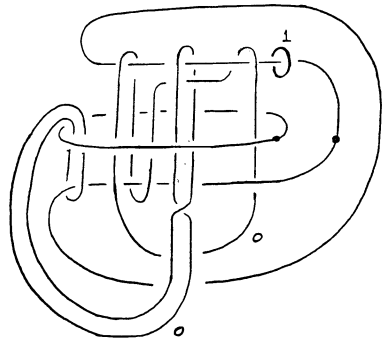
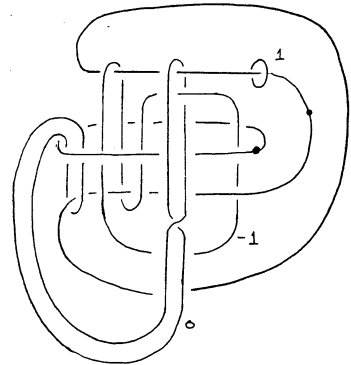

Figure 24.            Figure 25.

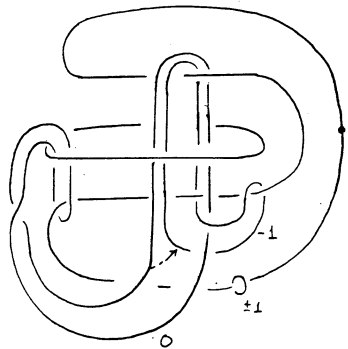
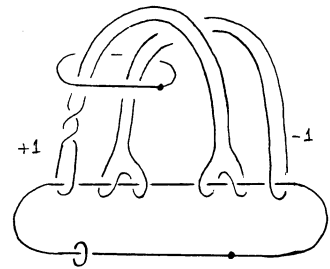

Figure 26.            Figure 27.

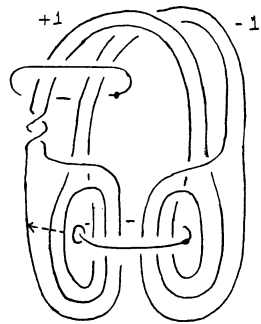

Figure 28.

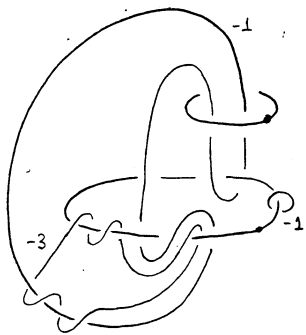
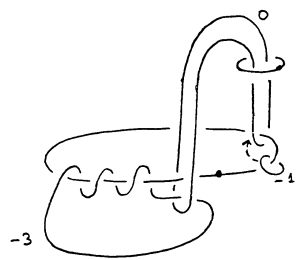

Figure 29.                                  Figure 30.

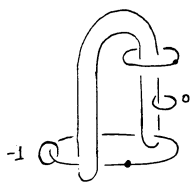
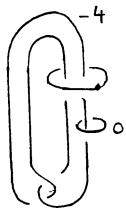
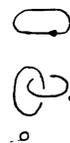

Figure 31.          Figure 32.          Figure 33.

*Remark* 1. It is remarkable to note that the above proof and the proof that the 2-fold covering $M$ of the Cappell-Shaneson fake $\mathbf{RP}^4$ ([CS2]) is diffeomorphic to $S^4$ ([G]) evolve similarly, i.e. by first showing that $M$ is obtained by Gluck twisting $S^4$, then canceling excess handles of $M$, and turning $M$ upside down ([AK2], [AK3]) (notice, after turning the handlebodies upside-down, the curious resemblence between Figure 27 of this paper and Figure 28 of [AK3]).

*Remark* 2. The complement of the imbedding $\Sigma \hookrightarrow S^3 \times S^1 \# S^2 \times S^2$ is given by a hyperbolic pair of imbedded 2-spheres, which can be described as $B^4$ with two 2-handles attached along two ribbon knots (Figure 8). However by [T] we know that we cannot find an imbedded $S^2$ repesenting one of the generators of this hyperbolic pair with simply connected complement.

*Remark* 3. Clearly there is more than one way to surger $\Sigma \times S^1$ to get $Q$; we chose to surger the loop in $\Sigma$ which corresponds the linking circle of the trefoil knot (Figure 1). In this paper we have not attempted to treat all the possibilities of surgering different loops as well as using different Rohlin invariant one homology spheres $\Sigma$, other then Poincaré homology sphere.

Also, there are two ways to surger a given loop in a 4-manifold corresponding to the parity of its framing; it can easily be checked that if we surger $\Sigma \times S^1$ along the same loop by using the other framing (i.e. by attaching a pair of linked 0-framed 2-handles in Figure 5, instead of $+1$ and 0 framed 2-handles) we get $S^3 \times S^1 \# \mathbf{CP^2} \# \mathbf{C\bar{P}^2}$.

Finally we would like use this opportunity to correct a minor mistake in the Figure of [Sc]: The 1 push-off of the trefoil knot is incorrectly drawn (to correct it, the two right-handed twists in the figure should be changed to two left-handed twists).

Michigan State University, East Lansing, MI
*E-mail address*: akbulut@math.msu.edu
M.S.R.I., Berkeley, CA
*E-mail address*: akbulut@msri.org

## References

[A1]  S. Akbulut, On 2-dimensional homology classes of 4-manifolds, Math. Proc. Cambridge Philos. Soc. **82** (1977), 99–106.

[A2]  ______, A fake 4-manifold, *Four-Manifold Theory* (Durham, NH, 1982), Contemp. Math. **35** (1984), 75–141.

[A3]  ______, On fake $S^3 \tilde{\times} S^1 \natural S^2 \times S^2$, *Combinatorial Methods in Topology and Algebraic Geom.* (Rochester, NY, 1982), Contemp. Math. **44** (1985), 281–286.

[A4]  ______, Constructing a fake 4-manifold by Gluck construction to a standard 4-manifold, Topology **27** (1988), 239–243.

[AK1] S. Akbulut and R. Kirby, Branched covers of surfaces in 4-manifolds, Math. Ann. **252** (1980), 111–131.

[AK2] ______, An exotic involution on $S^4$, Topology **18** (1979), 75–81.

[AK3] ______, A potential smooth counterexample in dimension 4 to the Poincaré conjecture, the Schoenflies conjecture, and the Andrew-Curtis conjecture, Topology **24** (1985), 375–390.

[CS1] S. Cappell and J. Shaneson, On four-dimensional surgery and applications, Comment. Math. Helv. **46** (1971), 500–528.

[CS2] ______, Some new four-manifolds, Ann. of Math. **104** (1976), 61–72.

[FP]  R. Fintushel and P. S. Pao, Identification of certain 4-manifolds with group actions, Proc. Amer. Math. Soc. **67** (1977), 344–350.

[G]   R. Gompf, Killing the Akbulut-Kirby 4-sphere, with relevance to the Andrews-Curtis and Schoenflies problems, Topology **30** (1991), 97–115.

[Sa]  Y. Sato, Scharlemann's 4-manifolds and smooth 2-knots in $S^2 \times S^2$, Proc. Amer. Math. Soc. **121** (1994), 1289–1294.

[Sc]  M. Scharlemann, Constructing strange manifolds with the dodecahedral space, Duke Math. J. **43** (1976), 33–40.

[T]   C. Taubes, Gauge theory on asymptotically periodic 4-manifolds, J. Differential Geom. **25** (1987), 363–430.